\documentclass{amsart}
\usepackage{amssymb}
\usepackage{amsfonts}

\newtheorem{theorem}{Theorem}
\theoremstyle{plain}

\newtheorem{corollary}{Corollary}

\newtheorem{definition}{Definition}

\newtheorem{lemma}{Lemma}

\numberwithin{equation}{section}
\input{tcilatex}

\begin{document}
\title[SOME HADAMARD-TYPE INEQUALITIES]{SOME HADAMARD-TYPE INEQUALITIES FOR
COORDINATED $P-$CONVEX FUNCTIONS AND GODUNOVA-LEVIN FUNCTIONS }
\author{$^{\blacklozenge ,\bigstar }$A. OCAK AKDEMIR}
\address{$^{\bigstar }$A\u{g}r\i\ \.{I}brahim \c{C}e\c{c}en University,
Faculty of Science and Arts, Department of Mathematics, 04100, A\u{g}r\i ,
Turkey}
\email{ahmetakdemir@agri.edu.tr}
\author{$^{\clubsuit }$M.EMIN OZDEMIR }
\curraddr{$^{\clubsuit }$Atat\"{u}rk University, K. K. Education Faculty,
Department of Mathematics, 25640, Campus, Turkey}
\email{emos@atauni.edu.tr}
\date{September, 22, 2010}
\keywords{Hadamard's inequality, co-ordinated convex, Godunova-Levin
functions, P- functions\\
$^{\blacklozenge }corresponding$ $author$}

\begin{abstract}
In this paper we established new Hadamard-type inequalities for functions
that co-ordinated Godunova-Levin functions and co-ordinated $P-$convex
functions, therefore we proved a new inequality involving product of convex
functions and $P-$functions on the co-ordinates$.$
\end{abstract}

\maketitle

\section{INTRODUCTION}

Let $f:I\subseteq 
\mathbb{R}
\rightarrow 
\mathbb{R}
$ be a convex function and let $a,b\in I,$ with $a<b.$ The following
inequality;%
\begin{equation*}
\ \ \ f(\frac{a+b}{2})\leq \frac{1}{b-a}\int_{a}^{b}f(x)dx\leq \frac{%
f(a)+f(b)}{2}
\end{equation*}

is known in the literature as Hadamard's inequality. Both inequalities hold
in the reversed direction if $f$ is concave.

In \cite{GL}, E.K. Godunova and V.I. Levin introduced the following class of
functions.

\begin{definition}
A function $f:I\subseteq 
\mathbb{R}
\rightarrow 
\mathbb{R}
$ is said to belong to the class of $Q(I)$ if it is nonnegative and, for all 
$x,y\in I$ and $\lambda \in (0,1)$ satisfies the inequality;%
\begin{equation*}
f(\lambda x+(1-\lambda )y)\leq \frac{f(x)}{\lambda }+\frac{f(y)}{1-\lambda }
\end{equation*}
\end{definition}

In \cite{DPP}, S.S. Dragomir et.al., defined following new class of
functions.

\begin{definition}
A function $f:I\subseteq 
\mathbb{R}
\rightarrow 
\mathbb{R}
$ is $P$ function or that $f$ belongs to the class of $P(I),$ if it is
nonnegative and for all $x,y\in I$ and $\lambda \in \lbrack 0,1],$ satisfies
the following inequality;

\begin{equation*}
f(\lambda x+(1-\lambda )y)\leq f(x)+f(y)
\end{equation*}
\end{definition}

In \cite{DPP}, S.S. Dragomir et.al., proved two inequalities of Hadamard's
type for class of Godunova-Levin functions and $P-$ functions.

\begin{theorem}
Let $f\in Q(I),a,b\in I,$ with $a<b$ and $f\in L_{1}[a,b].$ Then the
following inequality holds.
\end{theorem}

\begin{equation}
\ \ \ f(\frac{a+b}{2})\leq \frac{4}{b-a}\int_{a}^{b}f(x)dx  \label{a.1.1}
\end{equation}

\begin{theorem}
Let $f\in P(I),a,b\in I,$ with $a<b$ and $f\in L_{1}[a,b].$ Then the
following inequality holds.%
\begin{equation}
f(\frac{a+b}{2})\leq \frac{2}{b-a}\int_{a}^{b}f(x)dx\leq 2[f(a)+f(b)]
\label{a.1.2}
\end{equation}
\end{theorem}

In \cite{TUNC}, Tun\c{c} proved following theorem which containing product
of convex functions and $P-$functions.

\begin{theorem}
Let $a,b\in \left[ 0,\infty \right) $, $a<b,$ $I=\left[ a,b\right] $ with $%
f,g:\left[ a,b\right] \rightarrow 
\mathbb{R}
$\ be functions $f,g$ and $fg$ are in $L_{1}\left( \left[ a,b\right] \right) 
$. If $f$ is convex and $g$ belongs to the class of $P\left( I\right) $\
then,%
\begin{equation}
\frac{1}{b-a}\int_{a}^{b}f\left( x\right) g\left( x\right) dx\leq \frac{%
M\left( a,b\right) +N\left( a,b\right) }{2}  \label{a.1.3}
\end{equation}%
\textit{where }$M\left( a,b\right) =f\left( a\right) g\left( a\right)
+f\left( b\right) g\left( b\right) $\textit{\ and }$N\left( a,b\right)
=f\left( a\right) g\left( b\right) +f\left( b\right) g\left( a\right) $%
\textit{.}
\end{theorem}

In \cite{D1}, S.S. Dragomir defined convexity on the co-ordinates, as
following;

\begin{definition}
Let us consider the bidimensional interval $\Delta =[a,b]\times \lbrack c,d]$
in $%
\mathbb{R}
^{2}$ with $a<b,$ $c<d.$ A function $f:\Delta \rightarrow 
\mathbb{R}
$ will be called convex on the co-ordinates if the partial mappings $%
f_{y}:[a,b]\rightarrow 
\mathbb{R}
,$ $f_{y}(u)=f(u,y)$ and $f_{x}:[c,d]\rightarrow 
\mathbb{R}
,$ $f_{x}(v)=f(x,v)$ are convex where defined for all $y\in \lbrack c,d]$
and $x\in \lbrack a,b].$ Recall that the mapping $f:\Delta \rightarrow 
\mathbb{R}
$ is convex on $\Delta $ if the following inequality holds, 
\begin{equation*}
f(\lambda x+(1-\lambda )z,\lambda y+(1-\lambda )w)\leq \lambda
f(x,y)+(1-\lambda )f(z,w)
\end{equation*}%
for all $(x,y),(z,w)\in \Delta $ and $\lambda \in \lbrack 0,1].$
\end{definition}

Every convex function is co-ordinated convex but the converse is not
generally true.

In \cite{D1}, S.S. Dragomir established the following inequalities of
Hadamard's type for co-ordinated convex functions on a rectangle from the
plane $%
\mathbb{R}
^{2}.$

\begin{theorem}
Suppose that $f:\Delta =[a,b]\times \lbrack c,d]\rightarrow 
\mathbb{R}
$ is convex on the co-ordinates on $\Delta $. Then one has the inequalities;%
\begin{eqnarray}
&&\ f(\frac{a+b}{2},\frac{c+d}{2})  \label{a.1.6} \\
&\leq &\frac{1}{2}\left[ \frac{1}{b-a}\int_{a}^{b}f(x,\frac{c+d}{2})dx+\frac{%
1}{d-c}\int_{c}^{d}f(\frac{a+b}{2},y)dy\right]  \notag \\
&\leq &\frac{1}{(b-a)(d-c)}\int_{a}^{b}\int_{c}^{d}f(x,y)dxdy  \notag \\
&\leq &\frac{1}{4}\left[ \frac{1}{(b-a)}\int_{a}^{b}f(x,c)dx+\frac{1}{(b-a)}%
\int_{a}^{b}f(x,d)dx\right.  \notag \\
&&\left. +\frac{1}{(d-c)}\int_{c}^{d}f(a,y)dy+\frac{1}{(d-c)}%
\int_{c}^{d}f(b,y)dy\right]  \notag \\
&\leq &\frac{f(a,c)+f(a,d)+f(b,c)+f(b,d)}{4}  \notag
\end{eqnarray}
\end{theorem}

For recent results which similar to above inequalities see \cite{LA}, \cite%
{LAT}, \cite{DAR1}, \cite{DAR2} and \cite{DAR3}.

In \cite{OZ}, M.E. Ozdemir et.al., established the following Hadamard's type
inequalities as above for co-ordinated $m-$convex and $(\alpha ,m)-$convex
functions.

\begin{theorem}
Suppose that $f:\Delta =[0,b]\times \lbrack 0,d]\rightarrow 
\mathbb{R}
$ is $m-$convex on the co-ordinates on $\Delta $. If $0\leq a<b<\infty $ and 
$0\leq c<d<\infty $ with $m\in (0,1],$ then one has the inequality;%
\begin{eqnarray}
&&\frac{1}{(b-a)(d-c)}\int_{a}^{b}\int_{c}^{d}f(x,y)dxdy  \label{a.1.7} \\
&\leq &\frac{1}{4(b-a)}\min \left\{ v_{1},v_{2}\right\} +\frac{1}{4(d-c)}%
\min \left\{ v_{3},v_{4}\right\}  \notag
\end{eqnarray}%
where%
\begin{eqnarray*}
v_{1} &=&\int_{a}^{b}f(x,c)dx+m\int_{a}^{b}f(x,\frac{d}{m})dx \\
v_{2} &=&\int_{a}^{b}f(x,d)dx+m\int_{a}^{b}f(x,\frac{c}{m})dx \\
v_{3} &=&\int_{c}^{d}f(a,y)dy+m\int_{c}^{d}f(\frac{b}{m},y)dy \\
v_{4} &=&\int_{c}^{d}f(b,y)dy+m\int_{c}^{d}f(\frac{a}{m},y)dy.
\end{eqnarray*}
\end{theorem}

\begin{theorem}
Suppose that $f:\Delta =[0,b]\times \lbrack 0,d]\rightarrow 
\mathbb{R}
$ is $m-$convex on the co-ordinates on $\Delta $. If $0\leq a<b<\infty $ and 
$0\leq c<d<\infty $ $,$ $m\in (0,1]$ with $f_{x}\in L_{1}[0,d]$ and $%
f_{y}\in L_{1}[0,b],$ then one has the inequalities;%
\begin{eqnarray}
&&\frac{1}{b-a}\int_{a}^{b}f(x,\frac{c+d}{2})dx+\frac{1}{d-c}\int_{c}^{d}f(%
\frac{a+b}{2},y)dy  \label{a.1.8} \\
&\leq &\frac{1}{(b-a)(d-c)}\left[ \int_{a}^{b}\int_{c}^{d}\frac{f(x,y)+mf(x,%
\frac{y}{m})}{2}dydx\right.  \notag \\
&&\left. +\int_{c}^{d}\int_{a}^{b}\frac{f(x,y)+mf(\frac{x}{m},y)}{2}dxdy%
\right]  \notag
\end{eqnarray}
\end{theorem}

Similar results can be found for $(\alpha ,m)-$convex functions in \cite{OZ}%
. In this paper we established new Hadamard-type inequalities for
Godunova-Levin functions and $P-$functions on the co-ordinates on a
rectangle from the plane $%
\mathbb{R}
^{2}$ and we proved a new inequality involving product of co-ordinated
convex functions and co-ordinated $P-$functions$.$

\section{MAIN\ RESULTS}

We define Godunova-Levin functions and $P-$functions on the co-ordinates as
the following:

\begin{definition}
Let us consider the bidimensional interval $\Delta =[a,b]\times \lbrack c,d]$
in $%
\mathbb{R}
^{2}$ with $a<b,$ $c<d.$ A function $f:\Delta \rightarrow 
\mathbb{R}
$ is said to belong to the class of $Q(I)$ if it is nonnegative and for all $%
(x,y),(z,w)\in \Delta $ and $\lambda \in (0,1)$ satisfies the following
inequality;%
\begin{equation*}
f(\lambda x+(1-\lambda )z,\lambda y+(1-\lambda )w)\leq \frac{f(x,y)}{\lambda 
}+\frac{f(z,w)}{1-\lambda }
\end{equation*}
\end{definition}

A function $f:\Delta \rightarrow 
\mathbb{R}
$ is said to belong to the class of $Q(I)$ on $\Delta $ is called
coordinated Godunova-Levin function if the partial mappings $%
f_{y}:[a,b]\rightarrow 
\mathbb{R}
,$ $f_{y}(u)=f(u,y)$ and $f_{x}:[c,d]\rightarrow 
\mathbb{R}
,$ $f_{x}(v)=f(x,v)$ are belong to the class of $Q(I)$ where defined for all 
$y\in \lbrack c,d]$ and $x\in \lbrack a,b].$

We denote this class of functions by $QX(f,\Delta ).$ If the inequality
reversed then $f$ is said to be concave on $\Delta $ and we denote this
class of functions by $QV(f,\Delta ).$

\begin{definition}
Let $f:\Delta =[a,b]\times \lbrack c,d]\rightarrow 
\mathbb{R}
$ be a $P-$function with $a<b,$ $c<d.$ If it is nonnegative and for all $%
(x,y),(z,w)\in \Delta $ and $\lambda \in (0,1)$ the following inequality
\end{definition}

holds:%
\begin{equation*}
f(\lambda x+(1-\lambda )z,\lambda y+(1-\lambda )w)\leq f(x,y)+f(z,w)
\end{equation*}

A function $f:\Delta \rightarrow 
\mathbb{R}
$ is said to belong to the class of $P(I)$ on $\Delta $ is called
coordinated $P-$function if the partial mappings $f_{y}:[a,b]\rightarrow 
\mathbb{R}
,$ $f_{y}(u)=f(u,y)$ and $f_{x}:[c,d]\rightarrow 
\mathbb{R}
,$ $\ f_{x}(v)=f(x,v)$ are $P-$functions where defined for all $y\in \lbrack
c,d]$ and $x\in \lbrack a,b].$

We denote this class of functions by $PX(f,\Delta ).$ We need following
lemma for our main theorem.

\begin{lemma}
Every $f$ function that belongs to the class $Q(I)$ is said to belongs to
class $QX(f,\Delta ).$
\end{lemma}

\begin{proof}
Suppose that $f:\Delta =[a,b]\times \lbrack c,d]\rightarrow 
\mathbb{R}
$ is said to belong to the class $Q(I)$ on $\Delta .$ Consider the function $%
f_{x}:[c,d]\rightarrow \lbrack 0,\infty ),$ $f_{x}(v)=f(x,v).$ Then $\lambda
\in (0,1)$ and $v_{1},v_{2}\in \lbrack c,d],$ one has: 
\begin{eqnarray*}
f_{x}(\lambda v_{1}+(1-\lambda )v_{2}) &=&f(x,\lambda v_{1}+(1-\lambda
)v_{2}) \\
&=&f(\lambda x+(1-\lambda )x,\lambda v_{1}+(1-\lambda )v_{2}) \\
&\leq &\frac{f(x,v_{1})}{\lambda }+\frac{f(x,v_{2})}{1-\lambda } \\
&=&\frac{f_{x}(v_{1})}{\lambda }+\frac{f_{x}(v_{2})}{1-\lambda }
\end{eqnarray*}%
which shows convexity of $f_{x}.$ The fact that $f_{y}:[a,b]\rightarrow 
\mathbb{R}
,$ $f_{y}(u)=f(u,y)$ is also convex on $[a,b]$ for all $y\in \lbrack c,d]$
goes likewise and we shall omit the details.
\end{proof}

The following inequalities is considered the Hadamard-type inequalities for
Godunova-Levin functions on the co-ordinates.

\begin{theorem}
Suppose that $f:\Delta =[a,b]\times \lbrack c,d]\rightarrow 
\mathbb{R}
$ is said to belong to the class $QX(f,\Delta )$ on the co-ordinates on $%
\Delta $ with $f_{x}\in L_{1}[c,d]$ and $f_{y}\in L_{1}[a,b],$ then one has
the inequalities:%
\begin{eqnarray}
&&\frac{1}{16}\ \left[ f(\frac{a+b}{2},\frac{c+d}{2})\right]  \label{a.2.3}
\\
&\leq &\frac{1}{8}\left[ \frac{1}{b-a}\int_{a}^{b}f(x,\frac{c+d}{2})dx+\frac{%
1}{d-c}\int_{c}^{d}f(\frac{a+b}{2},y)dy\right]  \notag \\
&\leq &\frac{1}{(b-a)(d-c)}\int_{a}^{b}\int_{c}^{d}f(x,y)dydx  \notag
\end{eqnarray}
\end{theorem}

\begin{proof}
Since $f:\Delta =[a,b]\times \lbrack c,d]\rightarrow 
\mathbb{R}
$ is said to belong to the class $QX(f,\Delta )$ on the co-ordinates it
follows that the mapping $g_{x}:[c,d]\rightarrow 
\mathbb{R}
,$ $g_{x}(y)=f(x,y)$ is Godunova-Levin function on $[c,d]$ for all $x\in
\lbrack a,b].$ Then by Hadamard's inequality (\ref{a.1.1}) one has:%
\begin{equation*}
g_{x}(\frac{c+d}{2})\leq \frac{4}{d-c}\int_{c}^{d}g_{x}(y)dy,\forall x\in
\lbrack a,b].
\end{equation*}%
\qquad That is,%
\begin{equation*}
f(x,\frac{c+d}{2})\leq \frac{4}{d-c}\int_{c}^{d}f(x,y)dy,\forall x\in
\lbrack a,b].
\end{equation*}%
Integrating this inequality on $[a,b],$ we have:%
\begin{equation}
\frac{1}{b-a}\int_{a}^{b}f(x,\frac{c+d}{2})dx\leq \frac{4}{(b-a)(d-c)}%
\int_{a}^{b}\int_{c}^{d}f(x,y)dydx  \label{a.2.4}
\end{equation}%
A similar argument applied for the mapping $g_{y}:[a,b]\rightarrow 
\mathbb{R}
,$ $g_{y}(x)=f(x,y),$ we get:%
\begin{equation}
\frac{1}{d-c}\int_{c}^{d}f(\frac{a+b}{2},y)dy\leq \frac{4}{(b-a)(d-c)}%
\int_{c}^{d}\int_{a}^{b}f(x,y)dxdy  \label{a.2.5}
\end{equation}%
Summing the inequalities (\ref{a.2.4}), and (\ref{a.2.5}) , we get the last
inequality in (\ref{a.2.3}).

Therefore, by Hadamard's inequality (\ref{a.1.1}) we also have:%
\begin{equation*}
f(\frac{a+b}{2},\frac{c+d}{2})\leq \frac{4}{d-c}\int_{c}^{d}f(\frac{a+b}{2}%
,y)dy
\end{equation*}%
and%
\begin{equation*}
f(\frac{a+b}{2},\frac{c+d}{2})\leq \frac{4}{b-a}\int_{a}^{b}f(x,\frac{c+d}{2}%
)dx
\end{equation*}%
which give, by addition the first inequality in (\ref{a.2.3}).

This completes the proof.
\end{proof}

\begin{corollary}
Suppose that $f:\Delta =[a,b]\times \lbrack a,b]\rightarrow 
\mathbb{R}
$ is said to belong to the class $QX(f,\Delta )$ on the co-ordinates$,$ then
one has the inequalities:%
\begin{eqnarray}
\frac{1}{16}\left[ f(\frac{a+b}{2},\frac{a+b}{2})\right] &\leq &\frac{1}{8}%
\left[ \frac{1}{b-a}\int_{a}^{b}\left\{ f(x,\frac{a+b}{2})+f(\frac{a+b}{2}%
,x)\right\} dx\right]  \label{a.2.6} \\
&\leq &\frac{1}{(b-a)^{2}}\int_{a}^{b}\int_{a}^{b}f(x,y)dydx  \notag
\end{eqnarray}
\end{corollary}

\begin{corollary}
In (\ref{a.2.3}), under the assumptions Theorem 4 with $\ f(x,y)=f(y,x)$ for
all $x\in \lbrack a,b]\times \lbrack a,b],$ we have:%
\begin{eqnarray*}
f(\frac{a+b}{2},\frac{a+b}{2}) &\leq &\frac{1}{4}\left[ \frac{1}{b-a}%
\int_{a}^{b}f(x,\frac{a+b}{2})dx\right] \\
&\leq &\frac{1}{(b-a)^{2}}\int_{a}^{b}\int_{a}^{b}f(x,y)dydx
\end{eqnarray*}
\end{corollary}

\begin{lemma}
Every $P-$functions are coordinated on $\Delta $ or belong to the class of $%
PX(f,\Delta ).$
\end{lemma}

\begin{proof}
Let $f$ be a $P-$function and defined by $f_{y}:[a,b]\rightarrow 
\mathbb{R}
,$ $f_{y}(u)=f(u,y)$ and $f_{x}:[c,d]\rightarrow 
\mathbb{R}
,$ $f_{x}(v)=f(x,v)$ where $y\in \lbrack c,d]$, $x\in \lbrack a,b]$ and $%
\lambda \in \lbrack 0,1],v_{1},v_{2}\in \lbrack a,b],$ then 
\begin{eqnarray*}
f_{x}(\lambda v_{1}+(1-\lambda )v_{2}) &=&f(x,\lambda v_{1}+(1-\lambda
)v_{2}) \\
&=&f(\lambda x+(1-\lambda )x,\lambda v_{1}+(1-\lambda )v_{2}) \\
&\leq &f(x,v_{1})+f(x,v_{2}) \\
&=&f_{x}(v_{1})+f_{x}(v_{2})
\end{eqnarray*}%
which shows convexity of $f_{x}.$ The fact that $f_{y}:[a,b]\rightarrow 
\mathbb{R}
,$ $f_{y}(u)=f(u,y)$ is also convex on $[a,b]$ for all $y\in \lbrack c,d]$
goes likewise and we shall omit the details.
\end{proof}

The following inequalities is considered the Hadamard-type inequalities for $%
P$- functions on the co-ordinates.

\begin{theorem}
Suppose that $f:\Delta =[a,b]\times \lbrack c,d]\rightarrow 
\mathbb{R}
$ is said to belong to the class $PX(f,\Delta )$ on the co-ordinates on $%
\Delta $ with $f_{x}\in L_{1}[c,d]$ and $f_{y}\in L_{1}[a,b],$ then one has
the inequalities:%
\begin{eqnarray}
f(\frac{a+b}{2},\frac{c+d}{2}) &\leq &\frac{1}{b-a}\int_{a}^{b}f(x,\frac{c+d%
}{2})dx+\frac{1}{d-c}\int_{c}^{d}f(\frac{a+b}{2},y)dy  \label{a.2.7} \\
&\leq &\frac{4}{(b-a)(d-c)}\int_{a}^{b}\int_{c}^{d}f(x,y)dydx  \notag \\
&\leq &\frac{2}{(b-a)}\left[ \int_{a}^{b}f(x,c)dx+\int_{a}^{b}f(x,d)dx\right]
\notag \\
&&+\frac{2}{(d-c)}\left[ \int_{c}^{d}f(a,y)dy+\int_{c}^{d}f(b,y)dy\right] 
\notag
\end{eqnarray}
\end{theorem}

\begin{proof}
Since $f:\Delta =[a,b]\times \lbrack c,d]\rightarrow 
\mathbb{R}
$ is said to belong to the class $PX(f,\Delta )$ on the co-ordinates it
follows that the mapping $g_{x}:[c,d]\rightarrow 
\mathbb{R}
,$ $g_{x}(y)=f(x,y)$ is $P-$ function on $[c,d]$ for all $x\in \lbrack a,b].$
Then by Hadamard's inequality (\ref{a.1.2}) one has:%
\begin{equation*}
f(x,\frac{c+d}{2})\leq \frac{2}{d-c}\int_{c}^{d}f(x,y)dy\leq 2\left[
f(x,c)+f(x,d)\right]
\end{equation*}%
\qquad Integrating this inequality on $[a,b],$ we have:%
\begin{eqnarray}
\frac{1}{b-a}\int_{a}^{b}f(x,\frac{c+d}{2})dx &\leq &\frac{2}{(b-a)(d-c)}%
\int_{a}^{b}\int_{c}^{d}f(x,y)dydx  \label{a.2.8} \\
&\leq &\frac{2}{b-a}\left[ \int_{a}^{b}f(x,c)dx+\int_{a}^{b}f(x,d)dx\right] 
\notag
\end{eqnarray}%
A similar argument applied for the mapping $g_{y}:[a,b]\rightarrow 
\mathbb{R}
,$ $g_{y}(x)=f(x,y),$ we get:%
\begin{eqnarray}
\frac{1}{d-c}\int_{c}^{d}f(\frac{a+b}{2},y)dy &\leq &\frac{2}{(b-a)(d-c)}%
\int_{c}^{d}\int_{a}^{b}f(x,y)dxdy  \label{a.2.9} \\
&\leq &\frac{2}{d-c}\left[ \int_{c}^{d}f(a,y)dy+\int_{c}^{d}f(b,y)dy\right] 
\notag
\end{eqnarray}%
Addition (\ref{a.2.8}) and (\ref{a.2.9}), we get:%
\begin{eqnarray*}
\frac{1}{(b-a)(d-c)}\int_{a}^{b}\int_{c}^{d}f(x,y)dydx &\leq &\frac{1}{2(b-a)%
}\left[ \int_{a}^{b}f(x,c)dx+\int_{a}^{b}f(x,d)dx\right] \\
&&+\frac{1}{2(d-c)}\left[ \int_{c}^{d}f(a,y)dy+\int_{c}^{d}f(b,y)dy\right]
\end{eqnarray*}%
Which gives the last inequality in (\ref{a.2.7}). We also have: 
\begin{subequations}
\begin{align*}
& \frac{1}{b-a}\int_{a}^{b}f(x,\frac{c+d}{2})dx+\frac{1}{d-c}\int_{c}^{d}f(%
\frac{a+b}{2},y)dy \\
& \leq \frac{4}{(b-a)(d-c)}\int_{c}^{d}\int_{a}^{b}f(x,y)dxdy
\end{align*}%
Which gives the mid inequality in (\ref{a.2.7}). By Hadamard's inequality we
also have: 
\end{subequations}
\begin{equation*}
f(\frac{a+b}{2},\frac{c+d}{2})\leq \frac{2}{b-a}\int_{a}^{b}f(x,\frac{c+d}{2}%
)dx
\end{equation*}%
and%
\begin{equation*}
f(\frac{a+b}{2},\frac{c+d}{2})\leq \frac{2}{d-c}\int_{c}^{d}f(\frac{a+b}{2}%
,y)dy
\end{equation*}%
Adding these inequalities we get,%
\begin{equation*}
f(\frac{a+b}{2},\frac{c+d}{2})\leq \frac{1}{b-a}\int_{a}^{b}f(x,\frac{c+d}{2}%
)dx+\frac{1}{d-c}\int_{c}^{d}f(\frac{a+b}{2},y)dy
\end{equation*}%
Which gives the first inequality in (\ref{a.2.7}). This completes the proof.
\end{proof}

\begin{theorem}
Let $a,b,c,d\in \left[ 0,\infty \right) $,$a<b$ and $c<d,$ $\Delta =\left[
a,b\right] \times \lbrack c,d]$ with $f,g:\Delta \rightarrow 
\mathbb{R}
$\ be functions $f,g$ and $fg$ are in $L_{1}\left( \left[ a,b\right] \times
\lbrack c,d]\right) $. If $f$ is co-ordinated convex and $g$ belongs to the
class of $PX\left( f,\Delta \right) ,$\ then one has the inequality;%
\begin{eqnarray}
&&\frac{1}{\left( d-c\right) \left( b-a\right) }\int_{a}^{b}\int_{c}^{d}f%
\left( x,y\right) g\left( x,y\right) dydx  \label{a.3.1} \\
&\leq &\frac{L(a,b,c,d)+M(a,b,c,d)+N(a,b,c,d)}{4}  \notag
\end{eqnarray}%
\textit{where}%
\begin{eqnarray*}
L(a,b,c,d) &=&f(a,c)g(a,c)+f(b,c)g(b,c)+f(a,d)g(a,d)+f(b,d)g(b,d) \\
M(a,b,c,d) &=&f(a,c)g(a,d)+f(a,d)g(a,c)+f(b,c)g(b,d)+f(b,d)g(b,c) \\
&&+f(b,c)g(a,c)+f(b,d)g(a,d)+f(a,c)g(b,c)+f(a,d)g(b,d) \\
N(a,b,c,d) &=&f(b,c)g(a,d)+f(b,d)g(a,c)+f(a,c)g(b,d)+f(a,d)g(b,c)
\end{eqnarray*}
\end{theorem}

\begin{proof}
Since $f$ is co-ordinated convex and $g$ belongs to the class of $PX\left(
f,\Delta \right) ,$ by using partial mappings and from inequality (\ref%
{a.1.3}), we can write%
\begin{equation*}
\frac{1}{d-c}\int_{c}^{d}f_{x}\left( y\right) g_{x}\left( y\right) dy\leq 
\frac{f_{x}\left( c\right) g_{x}\left( c\right) +f_{x}\left( d\right)
g_{x}\left( d\right) +f_{x}\left( c\right) g_{x}\left( d\right) +f_{x}\left(
d\right) g_{x}\left( c\right) }{2}
\end{equation*}%
That is%
\begin{equation*}
\frac{1}{d-c}\int_{c}^{d}f\left( x,y\right) g\left( x,y\right) dy\leq \frac{%
f\left( x,c\right) g\left( x,c\right) +f\left( x,d\right) g(x,d)+f\left(
x,c\right) g\left( x,d\right) +f\left( x,d\right) g\left( x,c\right) }{2}
\end{equation*}%
Dividing both sides of this inequality $(b-a)$ and integrating over $[a,b]$
respect to $x,$ we have%
\begin{eqnarray}
&&\frac{1}{\left( d-c\right) \left( b-a\right) }\int_{a}^{b}\int_{c}^{d}f%
\left( x,y\right) g\left( x,y\right) dydx  \label{a.3.2} \\
&\leq &\frac{1}{2\left( b-a\right) }\int_{a}^{b}f\left( x,c\right) g\left(
x,c\right) +\frac{1}{2\left( b-a\right) }\int_{a}^{b}f\left( x,d\right)
g(x,d)  \notag \\
&&+\frac{1}{2\left( b-a\right) }\int_{a}^{b}f\left( x,c\right) g\left(
x,d\right) +\frac{1}{2\left( b-a\right) }\int_{a}^{b}f\left( x,d\right)
g\left( x,c\right)  \notag
\end{eqnarray}%
By applying (\ref{a.1.3}) to each integral on right hand side of (\ref{a.3.2}%
) and using these inequalities in (\ref{a.3.2}), we get the required result
as following%
\begin{eqnarray*}
&&\frac{1}{\left( d-c\right) \left( b-a\right) }\int_{a}^{b}\int_{c}^{d}f%
\left( x,y\right) g\left( x,y\right) dydx \\
&\leq &\frac{f\left( a,c\right) g\left( a,c\right) +f\left( b,c\right)
g(b,c)+f\left( a,c\right) g\left( b,c\right) +f\left( b,c\right) g\left(
a,c\right) }{4} \\
&&+\frac{f\left( a,d\right) g\left( a,d\right) +f\left( b,d\right)
g(b,d)+f\left( a,d\right) g\left( b,d\right) +f\left( b,d\right) g\left(
a,d\right) }{4} \\
&&\frac{f\left( a,c\right) g\left( a,d\right) +f\left( b,c\right)
g(b,d)+f\left( a,c\right) g\left( b,d\right) +f\left( b,c\right) g\left(
a,d\right) }{4} \\
&&+\frac{f\left( a,d\right) g\left( a,c\right) +f\left( b,d\right)
g(b,c)+f\left( a,d\right) g\left( b,c\right) +f\left( b,d\right) g\left(
a,c\right) }{4}
\end{eqnarray*}%
By a similar argument, if we apply (\ref{a.1.3}) for $f_{y}(x)g_{y}(x)$ on $%
[a,b],$ we get the same result.
\end{proof}


\begin{thebibliography}{99}
\bibitem{GL} E.K. Godunova and V.I. Levin,Neravenstra dlja funccii sirokogo
klassa soderzascego vypuklye, monotonnye i nekotorye drugie vidy
funkaii,Vycislitel Mat. i Mt.

Fiz.,Mezvuzov Sb. Nauc. Trudov. MPGI, Moscow, \ 1985,138-142.

\bibitem{DPP} S.S. Dragomir, J. Pecaric and L.E. Persson, Some inequalities
of Hadamard Type, Soochow Journal of Mathematics, Vol.21, No:3, pp. 335-341,
July 1995.

\bibitem{D1} S.S. Dragomir ,On Hadamard's inequality for convex functions on
the coordinates in a rectangle from the plane, Taiwanese Journal of
Mathematics, 5(2001), 775-788.

\bibitem{OZ} M.E. Ozdemir, E. Set, M. Z. Sar\i kaya, Some new Hadamard's
type inequalities for co-ordinated $m-$convex and $(\alpha ,m)-$convex
functions, Submitted.

\bibitem{LA} M. A. Latif, M. Alomari, On Hadamard-type inequalities for $h-$%
convex functions on the co-ordinates, International Journal of Math.
Analysis, 3 (2009), no. 33, 1645-1656.

\bibitem{LAT} M. A. Latif, M. Alomari, Hadamard-type inequalities for
product two convex functions on the co-ordinates, International Mathematical
Forum, 4 (2009), no. 47, 2327-2338.

\bibitem{DAR1} M. Alomari, M. Darus, Hadamard-type inequalities for $s-$%
convex functions, International Mathematical Forum, 3 (2008), no. 40,
1965-1975.

\bibitem{DAR2} M. Alomari, M. Darus, Co-ordinated $s-$convex function in the
first sense with some Hadamard-type inequalities, Int. Journal Contemp.
Math. Sciences, 3 (2008), no. 32, 1557-1567.

\bibitem{DAR3} M. Alomari, M. Darus, The Hadamard's inequality for $s-$%
convex function of $2-$variables on the co-ordinates, International Journal
of Math. Analysis, 2 (2008), no. 13, 629-638.

\bibitem{TUNC} M. Tun\c{c}, On Some Hadamard Type Inequalities for Product
of Different Kinds of Convex Functions, RGMIA, Res. Rep. Coll., 13 (2010),
Article 5, ONLINE: http://ajmaa.org/RGMIA/v13n1.php
\end{thebibliography}
\end{document}